\documentclass[12pt]{article}

\usepackage{graphics}
\usepackage{amsmath}
\usepackage{amsfonts}
\usepackage{amsthm}

\newtheorem{theorem}{Theorem}
\newtheorem{proposition}[theorem]{Proposition}
\newtheorem{lemma}[theorem]{Lemma}

\newtheorem{corollary}[theorem]{Corollary}

\def\Prob{{\rm Prob}}
\def\Probcr{{\rm Prob}_{cr}}
\def\Tone{{\mathbb T}}
\def\Hone{{\mathbb H}}

\def\Ecr{E_{\rm cr}}

\def\cnctd{\leftrightarrow}

\def\es{\emptyset}

\def\pr{{\prime}}
\def\ppr{{\prime\prime}}

\def\CMPsh{Commun. Math. Phys.}

\def\ZWVGsh{Z.~Wahrsheinlichkeitstheor.~Verw.~Geb.}

\newcommand{\eqn}[2]{\begin{equation}\label{#1}#2\end{equation}}
\newcommand{\eqnst}[1]{\begin{equation*}#1\end{equation*}}
\newcommand{\eqnspl}[2]{\begin{equation}\begin{split}\label{#1}%
	#2\end{split}\end{equation}}
\newcommand{\eqnsplst}[1]{\begin{equation*}\begin{split}%
	#1\end{split}\end{equation*}}

\title{The lowest crossing in 2D critical percolation}

\author{J.~van den Berg\\CWI, Amsterdam
\and A.~A.~J\'arai%
\thanks{Supported in part by NSERC of Canada, and the Pacific
Institute for the Mathematical Sciences.}
\\University of British Columbia, Vancouver}

\date{}

\begin{document}

\maketitle

\begin{abstract}
We study the following problem for critical site percolation on
the triangular lattice. Let $A$ and $B$ be sites  on a horizontal
line $e$ separated by distance $n$. Consider, in the half-plane above
$e$, the lowest occupied crossing $R$ from the
half-line left of $A$ to the half-line right of $B$.
We show that the probability that $R$ has a site at
distance smaller than $m$ from $AB$ is of order
$(\log (n/m))^{-1}$, uniformly in $1 \leq m \leq n/2$.
Much of our analysis can be carried out for other two-dimensional
lattices as well.

\end{abstract}

\section{Introduction}
The idea of the ``lowest'' crossing between two boundary pieces
of a domain is a well known and useful tool in the study of
two-dimensional percolation. Here we are interested in the question
how close the lowest crossing comes to the intermediate boundary
piece it has to cross.
To be specific, we fix the domain to be a
half plane and the two boundary pieces to be two disjoint half lines.

\subsection{Statement of the main result}
Let $\Tone$ denote the triangular lattice.
We note that much of our discussion applies to other lattices as well.
We consider $\Tone$ as a subset of the Euclidian plane, in such a way
that the distance between two neighbour vertices
of $\Tone$ is 1, and the integer points on the $X$-axis $e$ 
are vertices of $\Tone$. For notational convenience we denote these vertices
on $e$ by $\ldots, -2, -1, 0, 1, 2, \ldots$.
Denote the site $0$ by $A$ and the site $n$
by $B$.
Let $\ell = (-\infty, A) \cap \Tone$, $r = (B, \infty) \cap \Tone$,
and let $\Hone$ be the half plane above (and including) $e$.
Each site
$v \in \Tone$ is \emph{occupied} with probability $p$ and
\emph{vacant} with probability $1 - p$, independently. The corresponding
probability measure is denoted by $\Prob_p$, and expectation by $E_p$.
If $S_1,S_2$ are sets of sites, we say that $S_1$ is connected to $S_2$,
or $S_1 \cnctd S_2$, if there is a path of occupied sites that starts
in $S_1$ and ends in $S_2$. We say that
$S_1 \cnctd S_2 \text{ inside $S_3$}$, if all sites of the path are
in $S_3$.

All constants below are strictly positive and finite.
We write $a_n \asymp b_n$ to denote that there are constants $C_1$
and $C_2$ such that $C_1 a_n \le b_n \le C_2 a_n$. The exact values of
constants denoted by $C_i$ are not important for us, and $C_i$ may have
a different value from place to place. \\
\emph{Remark:}In the remainder of this paper  `path' will always
mean `self-avoiding path' (that is, a path which does not visit the same
site more than once).

\smallskip\noindent
\emph{The lowest crossing.} Consider all occupied paths between
$\ell$ and $r$ that stay inside $\Hone$. If there is such a path,
then there is a unique one closest to $AB$, call it $R$.
(See~\cite[p.~317]{Grimmett2} and~\cite{Kesten book} for a
discussion of the lowest crossing.) If $R$ contains a site on $AB$,
we call it a \emph{contact point}.


\medskip\noindent
We are only interested in contact points at criticality.
This is
because for $p < p_c$ the probability of an occupied crossing from
$\ell$ to $r$ decays exponentially as $n \to \infty$. Also, it is
not hard to see that for $p > p_c$ the fraction of those points on
$AB$ which are contact points is typically bounded away from $0$.
From now on we set $p = 1/2$, the critical
probability for site percolation on $\Tone$. We write $\Probcr$ for
$\Prob_{1/2}$. We note that by a Russo-Seymour-Welsh (RSW)
argument~\cite[Section 11.7]{Grimmett2},~\cite[Theorem 6.1]{Kesten book},%
~\cite{Russo1,Russo2,Seymour-Welsh},
we have $\Probcr(\text{$R$ exists}) = 1$.

Our main result is the following theorem:

\begin{theorem}\label{thm:main theorem}
We have, uniformly in $1 \leq m \leq n/2$,

$$\Probcr(\text{\rm $R$ has distance $< m$ from $AB$})
\asymp (\log (n/m))^{-1}.
$$
\end{theorem}

This theorem immediately implies (take $m=1$) the following Corollary:

\begin{corollary}\label{corol}

$$\Probcr(\text{\rm $R$ has a contact point})
\asymp (\log n)^{-1}, \,\, n \geq 1.
$$
\end{corollary}

\smallskip\noindent
{\emph Remarks} \\
(i) We like to note here that it is not even a priori
obvious that this probability goes to $0$ as $n$ goes to $\infty$.

\smallskip\noindent
(ii)
The only prerequisites needed in the proof are classical percolation
results: the RSW techniques and the fact that $p_c = 1/2$.
We do not use SLE processes, which were introduced by Schramm and which
have, by the work of him and other mathematicians,  recently
led to enormous progress (see ~\cite{Smirnov-Werner} and the references
given there). In fact we hope that
Theorem 1 will be useful in the study of $\text{SLE}_6$.
To illustrate this, note that
Theorem~\ref{thm:main theorem} indicates that in the scaling limit
when the lattice spacing goes to $0$ and the length of $AB$ is kept
fixed (say 1), the distribution of the distance of $R$ from $AB$ satisfies
\[ \Probcr(R \text{ has distance $<a$ from $AB$})
     \asymp (\log(1/a))^{-1}, \,\, a < 1/2. \]
In the scaling limit $R$ corresponds with the boundary of the hull of the 
chordal $\text{SLE}_6$ process in the half-plane started from $0$
and stopped at the first time it hits $(1,\infty)$ (see ~\cite{Smirnov},
Corollary 5). In this way one should obtain an analog of Theorem 1 in
terms of $\text{SLE}_6$.
The existence of a direct proof for $\text{SLE}_6$ of such a result is
not known to us.
Apart from these considerations, we think that Theorem 1 is interesting in
itself.

\subsection{Notation, definitions and key ingredients}
The theorem follows from the proposition below. This proposition uses
the knowledge of the critical exponent
describing the scaling of the probability that there are two disjoint
occupied paths in $\Hone$ that start at $0$ and end at distance $n$.
First we give some more definitions and notation.

For $n \ge 1$ and $v \in AB$ define the set
\[ H_n(v) = \{ u \in \Hone : |u - v| < n \}, \]
where $|\cdot|$ is the graph distance from the origin. We are also
going to need the half-annulus
\[ H_{n,m}(v) \stackrel{\text{def}}{=} H_n(v) \setminus H_m(v)
   = \{ u \in \Hone : m \le |u - v| < n \}. \]
If $S$ is a set of sites we set
\[ \partial S = \text{the set of sites in $S$ that have a
   neighbour in $S^c \cap \Hone$}, \]
and
\[ \bar\partial S = \text{the set of sites in $S^c \cap \Hone$
   that have a neighbour in $S$}. \]
We define the event
\[ D_n(v) = \{ \text{$\exists$ two disjoint occupied paths from
    $\bar\partial \{ v \}$ to $\partial H_n(v)$} \}. \]
We set
\[ \rho(n) = \Probcr(D_n(0)). \]
It is clear that this quantity will be important  in our analysis:
for a site $v \in AB$ to be a contact point, there must be two disjoint
occupied paths from $\bar\partial \{ v \}$ to the sets $\ell$ and $r$
respectively; when $v$ is in the bulk of $AB$ both sets have distance
of order $n$ from $v$. \\
We also need a version of $D_n$ for
$H_{n,m}(v)$. For $1 \le m < n$ let
\eqnsplst
{ D_{n,m}(v) &= \{ \text{$\exists$ two disjoint occupied paths from
   $\bar\partial H_m(v)$ to $\partial H_n(v)$} \}, \\
  \rho(n,m) &= \Probcr(D_{n,m}(0)). }
We are going to need the
following lemma about $\rho$.

\begin{lemma}\label{lem:rho}
We have
\begin{enumerate}
\item[(i)] $\rho(n) \asymp n^{-1}, \,\, n > 1$,
\item[(ii)] $\rho(n,m) \asymp (n/m)^{-1}$ uniformly
  in $1 \le m < n$.
\end{enumerate}
\end{lemma}

\medskip\noindent
Finally we state the following proposition. First, let
\[ X_{n,m} = |\{ 0 \le k \le n/m : \text{$H_m(km)$ is visited by
   $R$} \}|, \quad 1 \le m \le n/2. \]

\begin{proposition}\label{prop}
Uniformly in $1 \le m \le n/2$, with $n$ a multiple of $m$, we have
\begin{itemize}
\item[(i)]  $\Ecr X_{n,m} \asymp 1$,
\item[(ii)] $\Ecr ( X_{n,m} \,|\, X_{n,m} \ge 1 ) \asymp \log (n/m)$,
\item[(iii)] $\Ecr X_{n,m}^2 \asymp \log (n/m)$,
\item[(iv)] $\Probcr(X_{n,m} \ge 1) \asymp (\log (n/m))^{-1}$,
\end{itemize}
\end{proposition}

\subsection{Outline}
The rest of the paper is organized as follows.
In Subsection 2.1 we prove  Lemma~\ref{lem:rho}.
In Subsection 2.2 we prove Proposition~\ref{prop} from which,
as we will see in Subsection 2.3, Theorem~\ref{thm:main theorem} follows
immediately.
The only part which
uses the lattice structure in an essential way is the proof of the lemma.
The rest can easily be modified to suit other 2-dimensional lattices.

\section{Proofs}
\label{sec:proof}

\subsection{Proof of Lemma~\ref{lem:rho}}
A slightly weaker form
of this Lemma is the special case $j=2$ in Theorem 3 of
a recent preprint by Smirnov and Werner (2001) who use the recently
developed SLE machinery to derive this and many other results.
The proof below gives the somewhat stronger form we need, and is 
self-contained. 
(Also note that, for this special case, this stronger form gives an
answer to question 3 in Section 5 of the above mentioned paper by Smirnov and
Werner).
It is based on ideas
from~\cite[Section 2]{Kesten scaling},~\cite[Lemma 5]{KSZ} and
the unpublished work~\cite{Zhang crossing}, but we can bypass the use
of an $(\eta,k)$-fence, a notion introduced
in~\cite[Lemma 4]{Kesten scaling}.

\medskip\noindent
For $-n/2 \le k \le n/2 - 1$ define the events
\eqnsplst
{ P_{k,n} &= \left\{ \parbox{7.8cm}
   {$\exists$ occupied path from $k$ to
   $(-\infty,-n)$ and vacant path from $k + 1$ to $(n,\infty)$
   inside $\Hone$} \right\}, \\
  Q_{k,n} &= \left\{ \parbox{8.1cm}
   {$\exists$ disjoint occupied paths from $k$ to
   $(-\infty,-n)$ and from $k + 1$ to $(n,\infty)$ inside $\Hone$}
\right\}. }
If the event $P_{k,n}$  (or $Q_{k,n}$) occurs, let $S_1$ denote
the occupied path from $k$ to $(-\infty,-n)$ closest to $-n$. We claim
that $\Probcr(P_{k,n}) = \Probcr(Q_{k,n})$. Condition on $S_1$ and the
configuration ``below'' it. Then, since $p_c = 1/2$, flipping the
rest of the configuration establishes a one-to-one measure-preserving
correspondence between the two events.

We call a path $\pi$ in the half-annulus $H_{n,m}(v)$ a
\emph{half-circuit}, if it connects the two boundary pieces of
$H_{n,m}(v)$ lying on the boundary of $\Hone$. Let
\[ F_{n,m}(v) = \{ \text{$\exists$ occupied half-circuit in
   $H_{n,m}(v)$} \}. \]

\noindent
Further, let $P_n = \cup_{-n/2 \le k \le n/2 - 1} P_{k,n}$. 
Suppose there exists an occupied and a vacant path from $[-n/2,n/2]$ to
$(-\infty,-n)$ and to $(n,\infty)$ respectively. By considering the highest
such paths it is not difficult to see that then $P_n$ holds. Similarly we
see that the $P_{k,n}, \, -n/2 \le k \le n/2 - 1$,  are disjoint. So we have

\eqnspl{eq:bounds on Pn}
{ 1 &\ge \Probcr(P_n) = \sum_{-n/2 \le k \le n/2 - 1} \Probcr(P_{k,n})
   = \sum_{-n/2 \le k \le n/2 - 1} \Probcr(Q_{k,n}) \\
   &\ge \Probcr( F_{n-1,n/2}(-n) \cap \{ \text{$\exists$ vacant half
   circuit in $H_{n-1,n/2}(n)$} \}) \\
   &\ge C_1 > 0. }
The second of the inequalities follows because the two events on the
right hand side of this inequality imply (by the argument preceding
~\eqref{eq:bounds on Pn})
that $P_n$ occurs. The third inequality follows by independence and
the RSW Lemma.

\noindent
Since $Q_{k,n}$ is clearly at most $\rho(n/2)$, we get from 
~\ref{eq:bounds on Pn} that

\eqn{eq:rho lower}
{ \rho(n/2) \geq C_1/n.} 

\noindent
Further it is easy to see that for each $k \in [-n/2,n/2)$,
$Q_{k,n}$ contains the event 

\[ D_{4 n}(k) \cap F_{4n,3n}(k)
\cap \{ \text{all neighbours of $k$ occupied} \}. \]

\noindent
By FKG  and RSW this gives $\Probcr(Q_{k,n}) \geq C_2 \, \rho(4n)$.
Hence, by 
~\eqref{eq:bounds on Pn}),

\eqn{eq:rho upper}
{\rho(4n) \leq 1/(C_2 \, n).}

\noindent
Now ~\eqref{eq:rho upper} and
~\eqref{eq:rho lower} give 

\[ 1/(C_2 \, n) \geq \rho(4n) \geq C_1 / (8n). \]

\noindent
This (with the monotonicity of $\rho(n)$) gives immediately part (i) of the
Lemma.

%

\noindent
Part (ii) now follows from part (i) by a standard argument.
First of all, by inclusion of events and independence,
\[ \rho(n) \le \rho(m)\, \rho(n,m). \] 
To get an inequality in the
reverse direction we first note that we may assume that $2m \le n$.
It is not difficult to see that
\[ D_n \supset D_{2m} \cap F_{3m/2,m}(0) \cap F_{2m,3m/2}(0)
   \cap D_{n,m}. \]
By RSW, the second and third event on the r.h.s. are bounded away from $0$,
hence
\[ \rho(n) \ge C_4\, \rho(2 m) \rho(n,m). \]
This inequality, its above mentioned analog in the other direction, and part
(i) of the Lemma immediately gives part (ii).

\subsection{Proof of Proposition~\ref{prop}}

Let $R$, $A$ and $B$ be as in Section 1, and let
$1 \le m \le n/2$ with $n$ a multiple of $m$.
Observe that for $km \in AB$ we have
\eqn{eq:sufficient for visit}
{ \text{$R$ visits $H_m(km)$} \qquad \text{if and only if} \qquad
  \parbox{4.4cm}{$\exists$ occupied path from $\ell$ to $r$
     that visits $H_m(km)$,} }
and define the events
\eqnsplst
{ A_k &= \{ \text{$\exists$ occupied path from $\ell$ to $r$ that
    visits $H_m(km)$} \} \\
  &= \{ \text{$R$ visits $H_m(km)$} \}, \qquad 0 \le k \le
    n/m.}
We can write
\[ X_{n,m} = \sum_{0 \le k \le n/m} I[A_k], \]
where $I[\cdot]$ denotes the indicator of an event.

\smallskip\noindent
Throughout the proof we will assume that $m \geq 2$. The proof for $m = 1$
is similar and, in part (ii), simpler.

\begin{figure}
\center{\scalebox{.5}{%
  \includegraphics{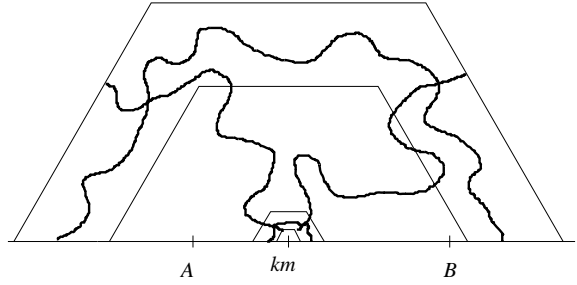}
}}
\caption{The events that force the occurrence of $A_k$.}
\label{fig1}
\end{figure}

\smallskip\noindent
\emph{Proof of (i).} We start with a lower bound for $\Ecr X_{n,m}$.
By 
inclusion of events (see Figure~\ref{fig1})
and the FKG inequality we have (with $F_{n,m}$ as in Section 2.1)
\eqnspl{eq:unconditional lower bound}
{ \Probcr(A_k) &\ge \Probcr(F_{2n,n}(km) \cap D_{2n,m/2}(km)
   \cap F_{m,m/2}(km)) \\
  &\ge \Probcr(F_{2n,n}(km))\, \rho(2n,m/2)\, \Probcr(F_{m,m/2}(km)). }
Here and later fractions are meant to be replaced by their integer
parts whenever necessary. By an RSW argument the first and third
factors are bounded below by some constant $C_1$. Therefore, by
Lemma~\ref{lem:rho} we have
\[ \Ecr X_{n,m} = \sum_{0 \le k \le n/m} \Probcr(A_k)
   \ge C_1^2 C_2 (n/m) (n/m)^{-1} = C_1^2 C_2. \]

For the upper bound we introduce the event
\[ G_{n,m}(v) = \{ \text{$\exists$ occupied path from
   $\bar\partial H_m(v)$ to $\partial H_n(v)$} \},
   \quad 1 \le m < n. \]
By an RSW argument
\eqn{eq:one path bound}
{ \Probcr(G_{n,m}) \le C_3 (n/m)^{-\mu} }
for some positive constants $\mu$ and $C_3$.
Let $1 \le k \le \frac{1}{2}(n/m)$,
and assume that the event $A_k$ occurs. Then it is easy to see that
the events $D_{km,m}(km)$
and $G_{n/2,km}(km)$ both occur. Since these latter events are
independent we have, by Lemma~\ref{lem:rho} and~\eqref{eq:one path bound},
\[ \Probcr(A_k) \le \Probcr(D_{km,m}(km))\, \Probcr(G_{n/2,km}(km))
   \le C_4 \frac{1}{k} \left( \frac{km}{n} \right)^{\mu}. \]
The sum of the right hand side over these $k$'s is bounded by some constant
$C_5$. A similar argument applies when $\frac{1}{2} (n/m) < k \le
(n/m) - 1$. Finally, in the cases $k = 0$ or $k =  n/m $
we have $\Probcr(A_k) \le 1$. This proves that $\Ecr X_{n,m} \le C_6$.

\smallskip\noindent
\emph{Proof of the lower bound in part (ii)}. \\
The idea in this proof is, roughly speaking, as follows: if $A_k$
occurs, there are from $H_m(k m)$ disjoint occupied paths to $\ell$ and $r$
respectively. Hence, to `let also $A_j$ occur' it (almost) suffices to
have two disjoint occupied paths from $H_m(j m)$ to the latter
path, and this should, by RSW arguments `cost' a probability of order
$\Probcr (D_{(j-k)m, m}(j m))$, which by the Lemma is of order $1/(j-k)$.
However, if one does the conditioning in a naive way, technical
difficulties arise because `negative information can seep through'.
Therefore the argument has to be done very carefully and an auxiliary
event (which we will call $F^*_k$ below) has to be introduced to
`neutralise' this negative information. We now give the precise arguments:

\medskip\noindent
Let $V$ denote the first intersection of $R$ with the set
\[ U = \bigcup_{0 \le k \le  n/m } H_m(km), \]
when $R$ is traversed from left to right. For $v \in \partial U$ let
$B_v = \{ V = v \}$, and define $k$ to be the index for which
$v \in H_m(km)$, choosing the smaller if there are two of them.
We prove the lower bound
\eqn{eq:lower-bnd}
{ \Probcr (A_j \,|\, B_v) \ge \frac{C_1}{j - k} \qquad
   \text{for $k + 4 \le j \le n/m - 1$,
   $1 \le k \le n/(2m)$.} }


Let
\eqnspl{eq:def-of-S_1}
{ R_1 &= \text{the piece of $R$ to the left of $V$, including
   the site $V$}, \\
  S_1(v) &= \parbox[t]{11cm}
   {lowest occupied path from $\ell$ to $v$ that is
   disjoint from $U$, apart from the site $v$.} }

We claim that on the event $B_v$ we have $R_1 = S_1(v)$. Since
$V = v$, we have that $R_1$ is disjoint from $U$, apart from $v$. If
$S_1(v)$ was lower than $R_1$, then we could use $S_1(v)$ and the piece
of $R$ to the right of $v$ to construct an occupied crossing lower
than $R$, a contradiction.

\smallskip\noindent
The proof of the lower bound in (ii) is based on the following observation.
\eqn{eq:decomp}
{ B_v = \bigcup_{\pi_1} \{ S_1(v) = \pi_1 \} \cap \Theta(\pi_1,v)
   \cap \Delta(\pi_1,v), }
where
\eqnsplst
{ \Theta(\pi_1,v) &= \left\{ \parbox{7.9cm}
    {$\exists$ vacant path $\pi_2^*$ from $\bar\partial\{v\}$ to
    $AB$, s.t.~$\pi_1$ is the occupied path from $\ell$ to $v$
    closest to $\pi_2^*$} \right\}, \\
  \Delta(\pi_1,v) &= \{\text{$\exists$ occupied path $\pi_3$ from
    $\bar\partial\{v\}$ to $r$ disjoint from $\pi_1$}\}, }
and where the union is over all paths $\pi_1$ from $\ell$ to
$v$ which are
disjoint from $U$, apart from the site $v$.
We will, for the time being, consider $v$ as fixed, and, to simplify
notation, write $S_1$, $\Theta(\pi_1)$ and $\Delta(\pi_1)$ instead of
$S_1(v)$ etc.

\smallskip\noindent
We first show that if $B_v$ occurs, then the right hand side
of~\eqref{eq:decomp} occurs. Take $\pi_1 = R_1$, then by the
discussion following~\eqref{eq:def-of-S_1} the event
$\{ S_1 = \pi_1 \}$ occurs. Since $R$ is the lowest crossing,
there is a vacant path from $\bar\partial\{v\}$ to $AB$. Take
$\pi_2^*$ to be the one closest to $\pi_1$. We claim that then also
$\pi_1$ is the occupied path closest to $\pi_2^*$. Let $\rho$ be an
occupied path from $\ell$ to $v$ that is closer to $\pi_2^*$ than
$\pi_1$. Since $\pi_2^*$ is below $R$, also $\rho$ is below $R$.
Now $\rho$ together with the piece of $R$ to the right of $v$ forms an
occupied crossing lower than $R$, a contradiction. This shows that
$\Theta(\pi_1)$ occurs. Finally, taking $\pi_3$ to be the piece of
$R$ to the right of $v$ shows that $\Delta(\pi_1)$ occurs.

Next assume that the right hand side of~\eqref{eq:decomp} occurs,
and choose the paths $\pi_1$, $\pi_2^*$ and $\pi_3$ that show this.
The fact that $\pi_1$, $\pi_3$ are occupied and that $\pi_2^*$ is
vacant implies that $R$ exists and passes through $v$. Thus $R_1$,
the piece of $R$ to the left of $v$, is defined. Also, $R$ lies below
the concatenation of $\pi_1$ and $\pi_3$. Since $\pi_2^*$ is vacant,
$R_1$ lies between $\pi_1$ and $\pi_2^*$. Since $\Theta(\pi_1)$ occurs,
$R_1 = \pi_1 = S_1$, and hence $v$ is the first intersection of $R$
with $U$, that is $B_v$ occurs.

Now we are ready to start the argument for~\eqref{eq:lower-bnd}.
By ~\eqref{eq:decomp} we can write
\eqn{eq:sum-pi1}
{ \Probcr (A_j \cap B_v)
  = \sum_{\pi_1} \Probcr(\{ S_1 = \pi_1 \} \cap \Theta(\pi_1)
   \cap \Delta(\pi_1) \cap A_j). }
Fix $\pi_1$, and
on the event
$\Delta(\pi_1)$
let $S_3(\pi_1)$ denote the highest occupied path from
$\bar\partial\{v\}$
to $r$ disjoint from $\pi_1$.
The occurrence of the event $\{ S_1 = \pi_1 \}$ only depends on the
states of $v$ and the sites that are on or below $\pi_1$ but outside $U$.
Let $\Omega(\pi_1)$
denote this set. For fixed $\pi_1$ the occurrence of
$\{ S_3(\pi_1) = \pi_3 \}$ only depends on sites above the union of
$\pi_1$ and $\pi_3$, and on the sites on $\pi_3$. Let
$\Omega(\pi_1,\pi_3)$ denote this set. (It may happen, but is not harmful,
that $\Omega(\pi_1) \cap \Omega(\pi_1,\pi_3) \not= \es$.)
We have
\[ \Delta(\pi_1) = \bigcup_{\pi_3} \{ S_3(\pi_1) = \pi_3 \}. \]
Thus we can write
\eqn{eq:sum-pi2}
{ \Probcr (A_j \cap B_v)
  = \sum_{\pi_1} \sum_{\pi_3}
   \Probcr(\{ S_1 = \pi_1,\, S_3(\pi_1) = \pi_3 \}
   \cap \Theta(\pi_1) \cap A_j). }

Now we construct events $K_{k,j}$ and $F^*_k$ such that the events
$K_{k,j}$ and $\{ S_1 = \pi_1,\, S_3(\pi_1) = \pi_3 \} \cap \Theta(\pi_1)$
are conditionally independent given $F^*_k$,
and moreover (on the event $B_v$)
$K_{k,j}$ forces the occurrence of $A_j$.
Let $\omega$ denote the configuration of occupied and vacant sites
in $\Hone$, and define the configuration $\omega^\pr$ by setting it
equal to a new independent configuration on
$\Omega(\pi_1) \cup \Omega(\pi_1,\pi_3)$, and equal to $\omega$
on $\Hone \setminus (\Omega(\pi_1) \cup \Omega(\pi_1,\pi_3))$. We let
\[ F^*_k = \{ \text{on $\omega^\pr$ $\exists$ vacant half-circuit
   in $H_{2m,m}(km)$} \}. \]
If $F^*_k$ occurs, then there is, in the configuration $\omega$,
a vacant path $\pi_4^*$
between $AB$ and $\pi_3$ creating a block. This means that
\eqn{eq:key-obs}
{ \parbox{10cm}{the path $\pi_2^*$ in the definition of
   $\Theta(\pi_1)$ can be chosen to lie on the left side
   of $\pi_4^*$.} }

\begin{figure}
\center{\scalebox{.5}{%
  \includegraphics{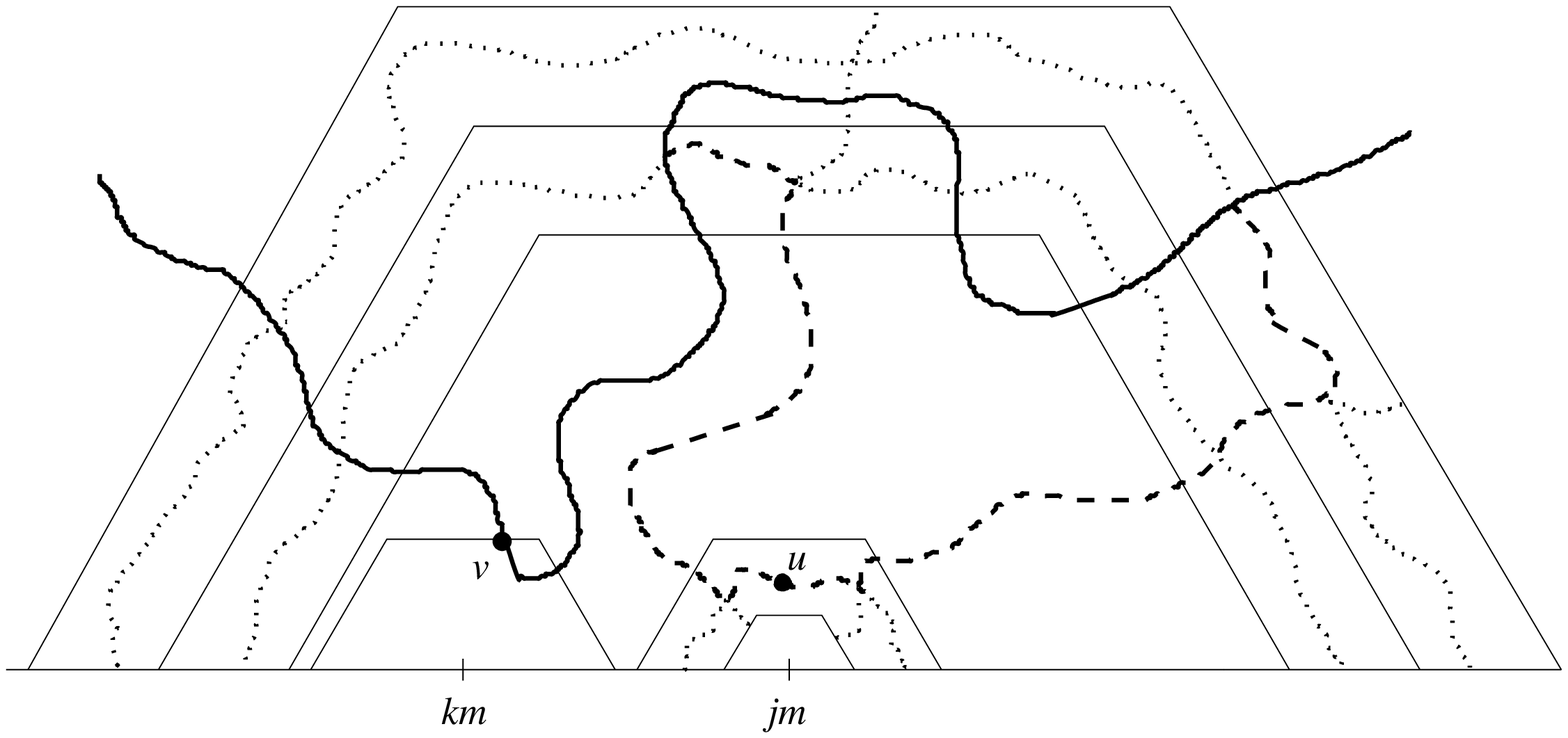}
}}
\caption{The dashed and dotted lines represent the event $K_{k,j}$
that forces the occurrence of $A_j$, given $B_v$. We used the dashed
parts to construct a path that visits $H_m(jm)$.}
\label{fig2}
\end{figure}

\smallskip\noindent
Next we define $K_{k,j}$
as the event that each of
the following four occurs on $\omega^\pr$:
\begin{itemize}
\item $\exists$ two disjoint occupied paths from
 $\bar\partial H_{m/2}(j m)$ to $\partial H_{4(j - k + 2)m}(j m)$
 that avoid the set $H_{2m}(km)$
\item $F_{4(j - k + 2)m, 2(j - k + 2)m}(j m)$
\item $F_{2(j - k + 2)m, (j - k + 2)m}(j m)$
\item $F_{m,m/2}(jm)$
\end{itemize}
We note that the first event we require is `almost'
$D_{4(j - k + 2)m, m/2}(jm)$. The only difference between these two events
is the avoidance condition, and it is easy to see that their
probabilities differ at most a constant factor.
Observe that if $K_{k,j}$ occurs, then
there is a path $\pi_5$ that is occupied on $\omega^\pr$, visits
$H_m(jm)$, and has both endpoints to the left of $H_m(km)$ on the
boundary of $\Hone$. Let $u$ be a site on $\pi_5$ that is in $H_m(jm)$.
If $u$ is above the union of $\pi_1$ and $\pi_3$ then $\pi_3$ visits
$H_m(jm)$. Otherwise there are points
$u^\pr, u^\ppr \in \pi_5 \cap \pi_3$ separated by $u$, which
implies that there is an occupied path (on $\omega$) from
$\bar\partial\{v\}$ to $r$ that visits $H_m(jm)$ (See Figure~\ref{fig2}).
Thus in both cases $A_j$ occurs.

By this observation and ~\eqref{eq:sum-pi2}, we have
\eqnspl{eq:extralabel}
{ \Probcr &(A_j \cap B_v) \\
  &\ge \sum_{\pi_1} \sum_{\pi_3}
   \Probcr(\{ S_1 = \pi_1,\, S_3(\pi_1) = \pi_3 \}
   \cap \Theta(\pi_1) \cap F^*_k \cap K_{k,j}). }
By~\eqref{eq:key-obs} and the construction of $K_{k,j}$ it
follows that, given  $F^*_k$, $K_{k,j}$ is conditionally
independent of
$\Theta(\pi_1) \cap \{ S_1 = \pi_1,\, S_3(\pi_1) = \pi_3 \}$.
Moreover, $K_{k,j}$ is independent of $F^*_k$.


This gives
that the right hand side of~\eqref{eq:extralabel} equals
\eqnspl{eq:use-indep}
{\sum_{\pi_1} \sum_{\pi_3}
   \Probcr(\{ S_1 = \pi_1,\, S_3(\pi_1) = \pi_3 \}
   \cap \Theta(\pi_1) \cap F^*_k)\, \Probcr(K_{k,j}). }
By the FKG inequality, Lemma 2 and RSW arguments
we have:
\eqn{eq:bnd-on-Kkj}
{ \Probcr(K_{k,j}) \ge C_2\, \rho(4(j-k+2)m,m) \ge \frac{C_3}{j-k}. }
To deal with the rest of the expression on the right hand side
of~\eqref{eq:use-indep} we condition on
the configuration $\sigma$ in $\Omega(\pi_1) \cup \Omega(\pi_1, \pi_3)$.
Note that, for fixed $\pi_1$,
$\pi_3$ and $\sigma$, the events $\Theta(\pi_1)$ and $F^*_k$ are
decreasing in the site variables in
$\Hone \setminus (\Omega(\pi_1) \cup \Omega(\pi_1,\pi_3))$.
Thus the FKG inequality implies that
\eqnspl{eq:bnd-FKG}
{ \Probcr& (\{ S_1 = \pi_1,\, S_3(\pi_1) = \pi_3 \}
   \cap \Theta(\pi_1) \cap F^*_k) \\
  &\ge \Probcr(\{ S_1 = \pi_1,\, S_3(\pi_1) = \pi_3 \}
   \cap \Theta(\pi_1)) \Probcr(F^*_k) \\
  &\ge C_4\, \Probcr(\{ S_1 = \pi_1,\, S_3(\pi_1) = \pi_3 \}
   \cap \Theta(\pi_1)). }
The bounds~\eqref{eq:extralabel},%
~\eqref{eq:use-indep},~\eqref{eq:bnd-on-Kkj}
and~\eqref{eq:bnd-FKG} (and~\eqref{eq:decomp})
yield
\eqnsplst
{ \Probcr& (A_j \cap B_v)
   \ge \frac{C_3\, C_4}{j - k} \sum_{\pi_1, \pi_3}
    \Probcr(\{ S_1 = \pi_1,\, S_3(\pi_1) = \pi_3 \}
    \cap \Theta(\pi_1)) \\
  &= \frac{C_3\, C_4}{j - k} \Probcr(B_v). }

Summing over $j$ this gives, for $v$ having $x$-coordinate
at most $n/2$,
\eqn{eq:cond lower}
{ \Ecr( X_{n,m} \,|\, B_v ) \ge C_3 \log (n/m). }
Let
\eqnst
{ J = \{ \text{$V$ has
$x-$coordinate $\le n/2$} \}
   = \bigcup_{v: v_x \le n/2} B_v,
}
where the union is over all $v \in \partial U$ with $x-$coordinate
at most $n/2$.
By symmetry, $\Probcr(J) \ge \frac{1}{2} \Probcr(X_{n,m} \ge 1)$. This
and~\eqref{eq:cond lower} gives
\eqnsplst
{ \Ecr& ( X_{n,m} \,|\, X_{n,m} \ge 1 )
    = \frac{\Ecr( X_{n,m} ; X_{n,m} \ge 1 )}{\Probcr(X_{n,m} \ge 1)}
    \ge \frac{\Ecr( X_{n,m} ; J )}{2\, \Probcr(J)} \\
   &= \frac{1}{2} \Ecr(X_{n,m} \,|\, J)
    \quad \ge (C_3/2) \log (n/m). }

\emph{Proof of the upper bound in (iii).}
In bounding $\Probcr(A_k \cap A_j)$ we may assume, by symmetry,
that $k \le j$ and $k \le n/m   - j $.
We may further assume that
$1 \le k \le j - 3$ by bounding $\Probcr(A_k \cap A_j)$ by
$\Probcr(A_j)$ in the cases $k = 0, j - 2, j - 1, j$ and using (i).
We separate three cases.

Case 1: $j - k < 2k$. Let $s = \lfloor (j - k - 1)/2 \rfloor$, and
$s^\pr = \lfloor (j - k)/2 \rfloor$. (We have $s^\pr = s$, if
$j - k$ is odd, and $s^\pr = s + 1$, if $j - k$ is even.) It is a
simple matter to check the inequalities $j - k \le k + s^\pr \le n/(2 m)$.
It is not difficult to see that if $A_k \cap A_j$ occurs,
then the following four events
occur:
\eqnst
{ D_{s m,m}(k m),\, D_{s m,m}(j m),\,
D_{(k + s^\pr) m,(j - k) m}((k + s^\pr) m),\,
   G_{n/2,(k + s^\pr) m}((k + s^\pr)m). }
Also note that these events are independent.
Thus by Lemma~\ref{lem:rho} and~\eqref{eq:one path bound}
\eqnspl{eq:case 1 bound}
{ \Probcr& (A_k \cap A_j) \le C_1 \frac{1}{s^2} \frac{j - k}{k + s^\pr}
   \left( \frac{(k + s^\pr)m}{n/2} \right)^\mu
  \le C_2 \frac{1}{(j - k)^2} \frac{j - k}{k}
   \left( \frac{k m}{n} \right)^\mu \\
  &= C_2\, (j - k)^{-1}\, k^{\mu - 1}\, \left( \frac{n}{m} \right)^{-\mu}, }
where at the second inequality we used $k \le k + s^\pr \le 2k$.
The sum of the right hand side of~\eqref{eq:case 1 bound} over $j$
is bounded by
$C_3\, (\log k)\, k^{\mu - 1}\, \left( \frac{n}{m} \right)^{-\mu}$. The sum
of this quantity over $k$ is bounded by
$C_4\, (\log (n/m)\, (n/m)^\mu\, (n/m)^{-\mu} = C_4 \log (n/m)$.

Case 2: $2k \le j - k \le 2(n/m - k)/3$. Define $s$ and $s^\pr$
as in Case 1. It is simple to check that $k \le s^\pr$ and
$k + s^\pr + (j - k) \le n/m$. In this case $A_k \cap A_j$ implies that
the following independent events occur:
\eqnst
{ D_{k m, m}(k m),\, G_{s^\pr m,k m}(k m),\, D_{s m, m}(j m),\,
G_{n - (k + s^\pr) m ,(j - k) m}((k + s^\pr) m). }
Thus we have
\eqnspl{eq:case 2 bound}
{ \Probcr& (A_k \cap A_j) \le C_5 \frac{1}{k}
   \left( \frac{k}{s^\pr} \right)^\mu \frac{1}{s}
   \left( \frac{j - k}{n/m - k - s^\pr} \right)^\mu \\
  &\le C_6 \frac{1}{k} \left( \frac{k}{j - k} \right)^\mu
   \frac{1}{j - k} \left( \frac{j - k}{n/m} \right)^\mu
  \le C_6\, k^{\mu - 1}\, (j - k)^{-1}\, (n/m)^{-\mu}, }
where in the second step we used that $n/m - k - s^\pr \ge n/(2 m)$.
The sum of the right hand side over $j$ is bounded by
$C_7\, (\log (n/m))\, k^{\mu - 1}\, (n/m)^{-\mu}$. The sum of this expression
over $k$ is bounded by $C_8\, (\log (n/m))\, (n/m)^\mu\, (n/m)^{-\mu}
= C_8\, \log (n/m)$.

Case 3: $j - k > 2(n/m - k)/3$. Our condition implies that (with $s$ and
$s'$ as before)
$k \le n/m - j < (j - k)/2$, hence $k \le n/m - j \le s$. This time
$A_k \cap A_j$ implies the
following independent events :
\eqnst
{ D_{k m, m}(k m),\, G_{s m,k m}(k m),\, D_{( n/m - j) m, m}(j m),
\, G_{s m,(n/m - j) m}(j m). }
This gives the bound
\eqnspl{eq:case 3 bound}
{ \Probcr& (A_k \cap A_j) \le C_9 \frac{1}{k}
   \left( \frac{k}{s} \right)^\mu \frac{1}{n/m - j}
   \left( \frac{n/m - j}{s} \right)^\mu \\
  &\le C_{10} \frac{1}{k}
   \left( \frac{k}{n/m} \right)^\mu \frac{1}{n/m - j}
   \left( \frac{n/m - j}{n/m} \right)^\mu \\
  &\le C_{10}\, k^{\mu - 1}\, (n/m - j)^{\mu - 1}\, (n/m)^{-2\mu}, }
where at the second inequality we used that
$s \ge (j - k - 2)/2 > (n/4 m) - 1$. The sum of the right hand side
of~\eqref{eq:case 3 bound} over $j$ and $k$ is bounded by some $C_{11}$.

The three cases and the remark about symmetry show that
\[ \Ecr X_{n,m}^2 =
    \sum_{0 \le j, k \le n/m} \Probcr(A_k \cap A_j)
    \le C_{12} \log (n/m). \]

\smallskip\noindent
\emph{Proof of (iv).}
From (i) and the lower bound in (ii) we get
\eqn{eq:upper bound on prob}
{ \Probcr(X_{n,m} \ge 1) =
\frac{\Ecr X_{n,m}}{\Ecr(X_{n,m} \,|\, X_{n,m} \ge 1)}
   \le \frac{C_1}{C_2 \log (n/m)}. }
On the other hand, by the Cauchy-Schwarz inequality
\eqn{eq:Cauchy-Schwarz}
{ \Ecr(X_{n,m}) = \Ecr(X_{n,m} I[X_{n,m} \ge 1])
   \le (\Ecr X_{n,m}^2)^{1/2} (\Probcr(X_{n,m} \ge 1))^{1/2}. }
The upper bound in (iii) and (i) imply
$\Probcr(X_{n,m} \ge 1) \ge C_3 (\log (n/m))^{-1}$.

\smallskip\noindent
\emph{Proof of the upper bound in (ii).}
The equality in~\eqref{eq:upper bound on prob} and (i) and (iv) now
give the upper bound in (ii).

\smallskip\noindent
\emph{Proof of the lower bound in (iii).}
Similarly,~\eqref{eq:Cauchy-Schwarz} and (i) and (iv) give the
lower bound in (iii). \qed


\subsection{Proof of Theorem 1}
The case where $n$ is a multiple of $m$ is
(by the definition of $X_{n,m}$)
clearly equivalent to
part (iv) of Proposition~\ref{prop}.
As to the general case, denote the probability in the statement of
the theorem by $f(n,m)$.
It is easy to see, using a simple RSW argument, that
if $n' < n < n' + m$,
then $f(n',m)$ and $f(n,m)$ differ at
most a factor $C>0$ which does not depend on $n$, $n'$ and $m$.
This observation, together with the special case, gives
the general case.

\bigbreak
\textbf{Acknowledgments.} We thank
Oded Schramm
for stimulating remarks that suggested to us that our
result could be interesting for SLE, and
B\'alint T\'oth
for a comment that led us to consider the present, quite general, form of
Theorem~\ref{thm:main theorem}. \\
A.J. also thanks CWI for its hospitality
during the summer of 2000.

\end{document}